\documentclass[amsthm]{elsart}

\usepackage{yjsco}
\usepackage{natbib}

\usepackage[latin1]{inputenc}
\usepackage{amssymb, amsmath}
\usepackage{epsfig, epic, eepic, url}

\begin{document}


\newcommand{\calA}{{\cal A}}
\newcommand{\calB}{{\cal B}}
\newcommand{\calF}{{\cal F}}
\newcommand{\calG}{{\cal G}}
\newcommand{\calR}{{\cal R}}
\newcommand{\calZ}{{\cal Z}}
\def\D{\mathcal{D}}
\def\L{\mathcal{L}}
\def\S{\mathcal{S}}
\def\I{\mathcal{I}}
\def\V{\mathcal{V}}
\def\E{\mathcal{E}}

\newcommand{\N}{{\mathbb N}}
\newcommand{\Z}{{\mathbb Z}}
\newcommand{\Q}{{\mathbb Q}}
\newcommand{\R}{{\mathbb R}}
\newcommand{\C}{{\mathbb C}}
\newcommand{\K}{{\mathbb K}}
\newcommand{\kk}{{\mathrm k}}

\def\J{\mathrm{J}}
\def\x{\mathrm{x}}
\def\a{\mathrm{a}}
\def\d{\mathrm{d}}

\def\GL{\mathrm{GL}}
\def\det{\mathrm{det}}
\def\SL{\mathrm{SL}}
\def\PSL{\mathrm{PSL}}
\def\PGL{\mathrm{PGL}}
\def\O{\mathrm{O}}

\def\gl{\mathfrak{gl}}
\def\g{\mathfrak{g}}
\def\h{\mathfrak{h}}

\newcommand{\Frac}[2]{\displaystyle \frac{#1}{#2}}
\newcommand{\Sum}[2]{\displaystyle{\sum_{#1}^{#2}}}
\newcommand{\Prod}[2]{\displaystyle{\prod_{#1}^{#2}}}
\newcommand{\Int}[2]{\displaystyle{\int_{#1}^{#2}}}
\newcommand{\Lim}[1]{\displaystyle{\lim_{#1}\ }}
\newcommand{\PD}[2]{\frac{\partial #1}{\partial #2}}
\newcommand{\pd}[1]{\frac{\partial}{\partial #1}}

\newenvironment{menumerate}{%
    \renewcommand{\theenumi}{\roman{enumi}}%
    \renewcommand{\labelenumi}{\rm(\theenumi)}%
    \begin{enumerate}} {\end{enumerate}}
     
\newenvironment{system}[1][]%
	{\begin{eqnarray} #1 \left\{ \begin{array}{lll}}%
	{\end{array} \right. \end{eqnarray}}

\newenvironment{meqnarray}%
	{\begin{eqnarray}  \begin{array}{rcl}}%
	{\end{array}  \end{eqnarray}}

\newenvironment{marray}%
	{\\ \begin{tabular}{ll}}
	{\end{tabular}\\}

\newenvironment{Pmatrix}
        {$ \left( \!\! \begin{array}{rr} } 
        {\end{array} \!\! \right) $}
   
\newtheorem{example}{Example}
\newtheorem{lemma}{Lemma}
\newtheorem{corollary}{Corollary}
\newtheorem{proposition}{Proposition}
\newtheorem{remark}{Remark}

\def\Hat{\widehat}
\def\Bar{\overline}

\def\fbar{{\bar f}}
\def\xbar{{\bar \x}}
\def\xb{{\bar x}}
\def\yb{{\bar y}}
\def\pb{{\bar p}}

\newcommand{\fleche}[3]{#1 \stackrel{#2}\longrightarrow #3}
\def\ssi{si et seulement si\ }
\newcommand{\tab}{\hspace*{\fill}}
\newcommand{\bs}{{\backslash}}
\newcommand{\eps}{{\varepsilon}}
\newcommand{\into}{{\;\rightarrow\;}}
\def\vect{\vec}
\newcommand{\afaire}{$$\vdots$$ \begin{center} {\sc a faire ...} \end{center} $$\vdots$$ }
\newcommand{\pref}[1]{(\ref{#1})}

\def\Maple{{\sc Maple}}
\def\RG{{\sc Rosenfeld-Gr\"obner}}

\newcommand{\algf}{\sffamily}
\newcommand{\BEGIN}{{\algf begin}}
\newcommand{\END}{{\algf end}}
\newcommand{\IF}{{\algf if}}
\newcommand{\THEN}{{\algf then}}
\newcommand{\ELSE}{{\algf else}}
\newcommand{\ELIF}{{\algf elif}}
\newcommand{\FI}{{\algf fi}}
\newcommand{\WHILE}{{\algf while}}
\newcommand{\FOR}{{\algf for}}
\newcommand{\DO}{{\algf do}}
\newcommand{\OD}{{\algf od}}
\newcommand{\RETURN}{{\algf return}}
\newcommand{\PROCEDURE}{{\algf procedure}}
\newcommand{\FUNCTION}{{\algf function}}
\newcommand{\INDENTER}{{\algf si} \=\+\kill}

\newcommand{\target}{\mathop{\mathrm{t}}}
\newcommand{\source}{\mathop{\mathrm{s}}}
\newcommand{\trdeg}{\mathop{\mathrm{tr~deg}}}
\newcommand{\jet}[2]{\jmath_{#1}^{#2}}
\newcommand{\rank}{\operatorname{rank}}
\newcommand{\sign}{\operatorname{sign}}
\newcommand{\ord}{\operatorname{ord}}
\newcommand{\aut}{\operatorname{aut}}
\newcommand{\Hom}{\operatorname{Hom}}
\newcommand{\codim}{\operatorname{codim}}
\newcommand{\coker}{\operatorname{coker}}
\newcommand{\rp}{\operatorname{rp}}
\newcommand{\leader}{\operatorname{ld}}
\newcommand{\card}{\operatorname{card}}
\newcommand{\Fr}{\operatorname{Frac}}
\newcommand{\RF}{\operatorname{\mathsf{reduced\_form}}}
\newcommand{\rang}{\operatorname{rang}}

\def \Id{\mathrm{Id}}

\def \diff{\mathrm{Diff}^{\mathrm{loc}} }
\def \diffg{\mathrm{Diff} }
\def \Esc{\mathrm{Esc}}

\newcommand{\initial}{\mathop{\mathsf{init}}}
\newcommand{\separant}{\mathop{\mathsf{sep}}}
\newcommand{\quo}{\mathop{\mathsf{quo}}}
\newcommand{\pquo}{\mathop{\mathsf{pquo}}}
\newcommand{\lcoeff}{\mathop{\mathsf{lcoeff}}}
\newcommand{\mvar}{\mathop{\mathsf{mvar}}}

\newcommand{\prem}{\mathop{\mathsf{prem}}}
\newcommand{\remp}{\mathrel{\mathsf{partial\_rem}}}
\newcommand{\remf}{\mathrel{\mathsf{full\_rem}}}
\renewcommand{\gcd}{\mathop{\mathrm{gcd}}}
\renewcommand{\deg}{\mathop{\mathrm{deg}}}
\newcommand{\pairs}{\mathop{\mathrm{pairs}}}
\newcommand{\dd}{\mathrm{d}}
\newcommand{\ideal}[1]{(#1)}
\newcommand{\cont}{\mathop{\mathrm{cont}}}
\newcommand{\pp}{\mathop{\mathrm{pp}}}
\newcommand{\pgcd}{\mathop{\mathrm{pgcd}}}
\newcommand{\ppmc}{\mathop{\mathrm{ppcm}}}
\newcommand{\init}{\mathop{\mathrm{initial}}}

\begin{frontmatter}
\title{New classification techniques for ordinary differential equations}

\thanks{This research was partly supported by ANR Gecko.}

\author{Raouf Dridi}

\address[2]{
 Département de Mathématiques,\\
Université Badji Mokhtar, \\
BP 12, Annaba, Algeria}
\ead{raouf.dridi@math.u-psud.fr}
\address[1]{Bureau 44, Bâtiment 425,\\
Département de Mathématiques,\\
91405 Orsay Cedex, France}

\author{Michel Petitot}
\address{Bureau 332, Bâtiment M3,\\
LIFL,  Université Lille I\\
59655 Villeneuve d'Ascq Cedex, France}
\ead{michel.petitot@lifl.fr}

\begin{abstract}
	The goal of the present paper is to propose an 
	enhanced ordinary differential equations solver by exploitation
	of the powerful equivalence method of \'Elie Cartan. 	 
	This solver returns a target equation equivalent to the equation to be solved and the  transformation realizing the equivalence.
	The target ODE is a member of  a dictionary of ODE, that are regarded as well-known, or at least well-studied.
	The dictionary considered in this article are ODE in a book of Kamke.  The major advantage of our solver is  that 
	the equivalence transformation	is obtained without integrating differential equations.
	We provide also a  theoretical contribution revealing
	the relationship between  the change of coordinates  that maps 
	two differential equations and their  symmetry pseudo-groups.\end{abstract}


\end{frontmatter}

\section{Introduction}
Current symbolic ODE solvers make use of a combination of
Lie symmetry methods and  pattern matching techniques.  
While  pattern matching techniques are used when the ODE matches
 a recognizable pattern (that is, for which a solving method is
already implemented), symmetry methods are reserved for
the non-classifiable cases \citep{Terrab97, Terrab98}. 
\begin{figure}[!hbt]\label{introfig} 
	\begin{center}
\setlength{\unitlength}{0.00083333in}
\begingroup\makeatletter\ifx\SetFigFont\undefined%
\gdef\SetFigFont#1#2#3#4#5{%
  \reset@font\fontsize{#1}{#2pt}%
  \fontfamily{#3}\fontseries{#4}\fontshape{#5}%
  \selectfont}%
\fi\endgroup%
{\renewcommand{\dashlinestretch}{30}
\begin{picture}(3851,3799)(0,-10)
\put(1875,885){\ellipse{2162}{656}}
\put(1787,2649){\ellipse{1324}{408}}
\path(1780,1599)(1780,1230)
\path(1765.275,1288.900)(1780.000,1230.000)(1794.725,1288.900)
\path(1817,539)(1817,256)(1144,256)
\path(1200.540,270.135)(1144.000,256.000)(1200.540,241.865)
\path(1817,256)(2558,256)
\path(2501.460,241.865)(2558.000,256.000)(2501.460,270.135)
\thicklines
\put(2603,161){\arc{90}{1.5708}{3.1416}}
\put(2603,352){\arc{90}{3.1416}{4.7124}}
\put(3220,352){\arc{90}{4.7124}{6.2832}}
\put(3220,161){\arc{90}{0}{1.5708}}
\path(2558,161)(2558,352)
\path(2603,397)(3220,397)
\path(3265,352)(3265,161)
\path(3220,116)(2603,116)
\put(796,3297){\arc{90}{1.5708}{3.1416}}
\put(796,3695){\arc{90}{3.1416}{4.7124}}
\put(2773,3695){\arc{90}{4.7124}{6.2832}}
\put(2773,3297){\arc{90}{0}{1.5708}}
\path(751,3297)(751,3695)
\path(796,3740)(2773,3740)
\path(2818,3695)(2818,3297)
\path(2773,3252)(796,3252)
\thinlines
\path(1780,3259)(1780,2843)
\path(1765.000,2901.500)(1780.000,2843.000)(1795.000,2901.500)
\thicklines
\put(104,104){\arc{120}{1.5708}{3.1416}}
\put(104,359){\arc{120}{3.1416}{4.7124}}
\put(1066,359){\arc{120}{4.7124}{6.2832}}
\put(1066,104){\arc{120}{0}{1.5708}}
\path(44,104)(44,359)
\path(104,419)(1066,419)
\path(1126,359)(1126,104)
\path(1066,44)(104,44)
\put(2896,2564){\arc{90}{1.5708}{3.1416}}
\put(2896,2820){\arc{90}{3.1416}{4.7124}}
\put(3762,2820){\arc{90}{4.7124}{6.2832}}
\put(3762,2564){\arc{90}{0}{1.5708}}
\path(2851,2564)(2851,2820)
\path(2896,2865)(3762,2865)
\path(3807,2820)(3807,2564)
\path(3762,2519)(2896,2519)
\thinlines
\path(2476,2669)(2815,2669)
\path(2768.500,2657.000)(2815.000,2669.000)(2768.500,2681.000)
\path(1801,2444)(1801,2039)
\path(1786.000,2084.000)(1801.000,2039.000)(1816.000,2084.000)
\thicklines
\put(1246,1664){\arc{90}{1.5708}{3.1416}}
\put(1246,1968){\arc{90}{3.1416}{4.7124}}
\put(2421,1968){\arc{90}{4.7124}{6.2832}}
\put(2421,1664){\arc{90}{0}{1.5708}}
\path(1201,1664)(1201,1968)
\path(1246,2013)(2421,2013)
\path(2466,1968)(2466,1664)
\path(2421,1619)(1246,1619)
\put(1604,3543){\makebox(0,0)[lb]{{\SetFigFont{8}{9.6}{\rmdefault}{\bfdefault}{\updefault}Input :}}}
\put(1958,2200){\makebox(0,0)[lb]{{\SetFigFont{7}{8.4}{\rmdefault}{\mddefault}{\itdefault}No}}}
\put(2806,220){\makebox(0,0)[lb]{{\SetFigFont{8}{9.6}{\rmdefault}{\bfdefault}{\updefault}Fail}}}
\put(2064,362){\makebox(0,0)[lb]{{\SetFigFont{7}{8.4}{\rmdefault}{\mddefault}{\itdefault}No}}}
\put(1497,362){\makebox(0,0)[lb]{{\SetFigFont{7}{8.4}{\rmdefault}{\mddefault}{\itdefault}Yes}}}
\put(2551,2744){\makebox(0,0)[lb]{{\SetFigFont{7}{8.4}{\rmdefault}{\mddefault}{\itdefault}Yes}}}
\put(1126,944){\makebox(0,0)[lb]{{\SetFigFont{8}{9.6}{\rmdefault}{\mddefault}{\itdefault} Does the equation admit}}}
\put(826,3344){\makebox(0,0)[lb]{{\SetFigFont{8}{9.6}{\rmdefault}{\mddefault}{\updefault} Ordinary differential equation  }}}
\put(226,119){\makebox(0,0)[lb]{{\SetFigFont{8}{9.6}{\rmdefault}{\bfdefault}{\updefault}of the order}}}
\put(226,269){\makebox(0,0)[lb]{{\SetFigFont{8}{9.6}{\rmdefault}{\bfdefault}{\updefault}Reduction }}}
\put(2926,2594){\makebox(0,0)[lb]{{\SetFigFont{8}{9.6}{\rmdefault}{\bfdefault}{\updefault}Integration  }}}
\put(1351,1694){\makebox(0,0)[lb]{{\SetFigFont{8}{9.6}{\rmdefault}{\bfdefault}{\updefault}Lie equations }}}
\put(1351,1844){\makebox(0,0)[lb]{{\SetFigFont{8}{9.6}{\rmdefault}{\bfdefault}{\updefault}Resolution of }}}
\put(1276,2594){\makebox(0,0)[lb]{{\SetFigFont{8}{9.6}{\rmdefault}{\mddefault}{\itdefault}Known equation?     }}}
\put(826,794){\makebox(0,0)[lb]{{\SetFigFont{8}{9.6}{\rmdefault}{\mddefault}{\itdefault} 1-parameter symmetry group?     }}}
\end{picture}
}
	\caption{General flowchart of typical ODE-solver.} \label{introfig0}
	\end{center}
\end{figure}
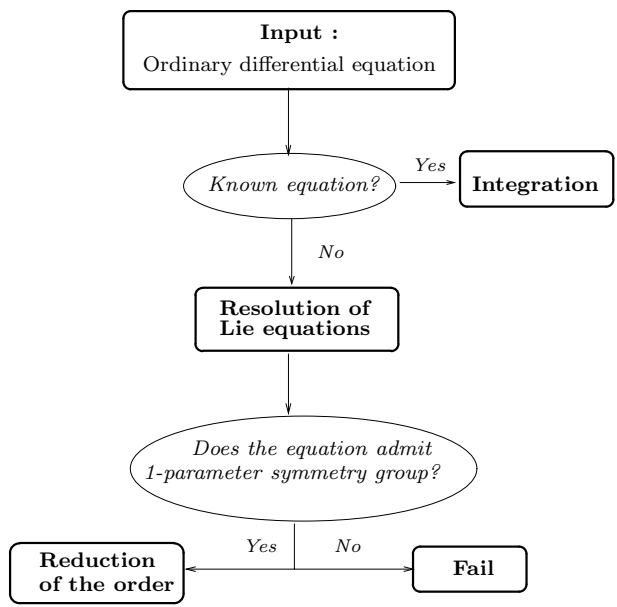 

Nevertheless, in practice,  these solvers often fail to return closed form solutions.  
This is the case for instance of  the equation
\begin{equation}\label{ODE1}
	 y''+ {y}^{3}{y'}^{4} + {\frac {{y'}^{2}}{y}} + \frac{1}{2}y =0,
\end{equation}
which   admits  only one 1-parameter symmetry  group. Using this information,
 present solvers return a complicated first order ODE
  and a quadrature which is quite useless for practical applications. 

More dramatically, when applied, to  the following  equation, these solvers output no result
\begin{equation}\label{ODE2}
	y'' - {\frac {2\,{x}^{4}y' - 6\,{y}^{2}x-1}{{x}^{5}}} =0.
\end{equation}
This failure is due to the fact that the above equation does not match any recognizable pattern
and has zero-dimensional point symmetry (pseudo-)group. Thus neither  symmetry  methods nor classification methods works. The same
goes for the equation
\begin{equation}\label{ODE3}
y''y + {y'}^{2}-4\,{y}^{6}+12\,x{y}^{4} - 
 \left( 4\,x + 12\,{x}^{2} \right) {y}^{2}+4\,{x}^{3}+4\,{x}^{2} -2\,\alpha= 0
\end{equation}

Our solver is designed to handle such differential equations.
It returns a target equation equivalent to the equation to be solved and the equivalence transformation. 
The target ODE is a member of  a dictionary of ODE, that are regarded as well-known, or at least well-studied.
The dictionary considered in the article are ODE in a book of Kamke.  However other ODE could be added 
to the dictionary without difficulty and this is an advantage of the method.

For the equation  \pref{ODE1}, we obtain the Rayleigh equation~(number 6.72, page 559 in {Kamke's} book)  
$$
	y'' + {y'}^4 + y=0
$$ 
and the change  of coordinates 
$$
	 (x,\, y) \into (x, \, \frac{y^2}{2}).
$$ 
For the equation \pref{ODE2}, we obtain the  first Painlev\'e equation (number 6.3, page 542 in {Kamke's} book) 
$$
	y'' =6y^2+x
$$ 
and the change of coordinates 
$$
	(x,\, y) \into (\frac{1}{x},\,  y).
$$
The equation \pref{ODE3} is mapped to  the  second Painlev\'e equation (number 6.6  in {Kamke's} enumeration) 
$$
	y''=2 y^{3}+yx+ \alpha
$$ 
under the transformation 
$$
	(x,\, y) \into (2x,\,  y^2 -x).
$$


To incorporate such changes of variables, one needs to  understand the \emph{equivalence problem} : Given two ODE  
 $$
 	E_f : \, y'' = f(x,y,y'), 
	\mbox{  } \, E_\fbar : \, \bar y''=\fbar (\bar x,\bar y,\bar y')
 $$
and an allowed  Lie pseudo-group of transformations $\Gamma\Phi$  acting on the variables $x$ and~$y$~(for real-life reasons, we restrict ourselves to second order ODE.
  However, what follows remains true for any order).  We shall say that $E_f$ and $E_\fbar$
are  equivalent under the action of the pseudo-group $\Gamma\Phi$  if there exists a change of coordinates $\varphi\in\Gamma\Phi$ that maps one equation to another. This  will be denoted by 
\begin{equation}\label{PDEsystem}
         E_\fbar = \varphi_*(E_f) \mbox{ and } \varphi\in \Gamma\Phi,
\end{equation}
or in abridged form $E_f\sim_\Phi E_\fbar$. As we shall see, the system \pref{PDEsystem} is PDE's system.
 One can consider the equivalence of an equation with its self : A \emph{symmetry}~$\sigma$ of the equation~$E_f$ is a solution of the self-equivalence obtained by setting $E_f=E_\fbar$ in the a PDE's system \pref{PDEsystem}. The  solutions of this self-equivalence problem form a Lie pseudo-group, the symmetry pseudo-group, which will be denoted by $S_{E_f,\Phi}$.

In practice, one distinguishes two possible situations in the computation of the change of coordinates.   First,  the input equation (the equation to be solved) and the target equation are known. This is an online computation. In the second situation, considered here in the construction of our solver, only the target equation
(an equation from Kamke's list) is known and we look for the change of coordinate that maps the generic equation to the target one.  So, assume that   $E_f$ is  a generic  differential equation  and $E_\fbar$ is fixed equation that
 falls within the effective differential algebra i.e. $E_\fbar$   is given by  equalities between differential 
polynomials with rational coefficients. 
   
Crucial in the construction of our solver, is the establishment of the relationship between the change of coordinates
 and  symmetry pseudo-groups. In particular, one might ask under which  conditions  the change of coordinates  
 can be obtained \emph{without} integrating differential equation ? To the authors's knowledge,   this is the
  first time in equivalence problem theory that such questions are investigated. 
  The answer, which constitutes the theoretical contribution of the paper, can be 
  be summarized as follows (we emphasis on  the fact that what follows remains true for any order) :

\begin{menumerate}

\item	[$(i)$]
The number of  constants appearing in the change of coordinates $\varphi\in\Gamma\Phi$,
mapping~$E_f$ to $E_\fbar$, is exactly the dimension of $S_{E_\fbar, \Phi}$.  
This implies that when this dimension vanishes the change of coordinate can be obtained without
  integrating differential equations.  Also, we have $\dim(S_{E_f, \Phi})=\dim(S_{E_\fbar, \Phi})$.

\item	[$(ii)$]
In the particular case when $\dim(S_{E_\fbar, \Phi})=0$,  the transformation $\varphi$ is algebraic in $f$
and its partial derivatives. The degree of this transformation $\varphi$ is exactly $\card(S_{E_\fbar, \Phi})$. 
In this case, the symmetry pseudo-groups $S_{E_f, \Phi}$ and $S_{ E_\fbar, \Phi}$ have finite cardinals. However,
 they need not to have the same cardinal.
\end{menumerate}

The simple fact  $\dim(S_{E_f, \Phi})=\dim(S_{E_\fbar, \Phi})$ allows us to construct a powerful hashing function
 which  significantly restricts the space of research in kamke's list. 
 For this reason, we use  7 possible types of transformations $\Phi_1,\cdots, \Phi_7$ (see table~{\ref{table}} page~\pageref{table}). 
We pre-calculate  to each target equation a \emph{signature index},
that is, the dimensions of the 7 symmetry pseudo-groups associated to the 7 types of transformations (this calculation
is done without any integration). If two differential equations are equivalent then their signature indices \emph{match}. 

The transformation $\varphi$ in  $(ii)$  can be obtained using differential elimination algorithms. 
This is explained in third section.  Unfortunately, such approach is rarely effective due to expressions swell. In order to avoid this,
 we propose in section~\ref{cartan} a new method to pre-compute the transformation $\varphi$ in terms of differential
invariants (we do this for each target equation $E_\fbar$ in Kamke's list). These invariants are provided by Cartan's  method.

\section{Equivalence problem}

The equivalence problem is the study of the action of a given pseudo--group of transformations on 
a given class of differential systems. In the algebraic framework, this action is viewed as the action of 
a $\D$--groupoid $\Phi$ on a diffiety $\E$.

\subsection{Equivalence problems and groupoids}
Recall that a diffiety  (see  \ref{annexe:diffalg}) is the set of formal Taylor series which are regular solutions
of a finite PDE's system. It is a pro-algebraic variety, fibered over an algebraic variety~$X$ and which will be denoted by 
$\pi: \E \longrightarrow X$. The projection of a Taylor series~${\jet \x \infty f \in \E}$, of a function $f: X \rightarrow U$, 
is the expansion point $x\in X$.  The   coordinate ring of a diffiety  is a reduced finitely generated \emph{differential} algebra. 
The automorphisms group of the diffiety is the set of the contact transformations (see \cite{Olver}) from $\E$ to $\E$.
 
A $\D$-groupoid $\Phi$ is a diffiety formed by invertible Taylor series and closed under the composition (see section  \ref{annexe:groupoids}).
A $\D$-groupoid acting on a  manifold $X$, is a subset of the space of infinite invertible jets $\J^\infty_* (X,X)$.
The Taylor series of contact transformations of a diffiety $\E$ form a $\D$-groupoid that acts on $\E$ and which will be denoted by~$\aut(\E)$.

Given a diffiety $\E$ fibred over $X$ and a $\D$-groupoid $\Phi$ acting on $X$. An equivalence problem is the action of $\Phi$ on the diffiety $\E$, 
that is, an injective representation\footnote{\'Elie Cartan call \emph{prolongement holohédrique} such injective representation.}, i.e. an injective morphism of $\D$-groupoids
$$
	\rho: \Phi \longrightarrow \aut(\E).
$$

\begin{example} [2nd order ODE, $\bar x = x + C, \bar y = \eta(x,y)$] \label{ex}

The infinite dimensional $\D$-groupoid $\Phi := \Phi_3$ acts on the points $(x,y) \in \J^0(\C, \C)$. 
The corresponding Lie defining equations are (we set $\xbar = \xi(x,y)$) 
\begin{equation} \label{Phi3:eqns}
       	    \xi_x = 1,\, \xi_y = 0,\, \eta_y \neq 0.
\end{equation} 
The action of $\Phi$ on the   jets space   $\J^0(\C, \C)$ is prolonged  (see appendix \ref{prolong}) to 
an  action on the first order jets space $X := \J^1(\C,\C)$. In the coordinates $(x,y,p=y')$, this action reads
\begin{equation} 
\begin{gathered}
    	\bar x = \xi,\ \bar y = \eta, \ 
		\bar p = \eta_x + p \eta_y.
\end{gathered}
\end{equation}
The equivalence condition \pref{PDEsystem} is obtained by prolonging the action of   $\Phi$
on the second order  jets space $\J^2(\C,\C)$ and by setting $y''=f(x,y,p)$. Thus, we obtain
\begin{equation} \label{ex:f}
\begin{gathered}
	    \xi_x = 1,\, \xi_y = 0,\, \xi_p=0, w\ \eta_y = 1, \\
	    \bar x = \xi,\ \bar y = \eta, \ \bar p = \eta_x + p \eta_y, \\
		\fbar(\xb, \yb, \pb) = \eta_{xx} + 2p \eta_{xy} + p^2 \eta_{yy} + f(x,y,p)\ \eta_y.
\end{gathered}
\end{equation}
This action of  $\Phi$ on the equation $E_f$ is viewed as an action on the Taylor series~${\jet \x \infty f}$, in other words,
as an action on the trivial diffiety  $\E := \J^\infty(X, \C)$ (fibred over the  manifold~$X$).
The coordinates ring of $X$ is  $\C[X] := \C[x,y,p]$ and the coordinate ring of~$\E$  is the ring 
of differential polynomials 
$$
	\C[\E] := \C[x,y,p]\{ f \} \mbox{ with } \Delta = \left\{ \pd{x}, \pd{y}, \pd{p} \right\}.
$$
\end{example}

\section{Differential-algebraic approach} \label{diffalg}
The aim of this section is  to use differential elimination  to solve the equivalence problem 
 \pref{PDEsystem} when the target function $\fbar(\xb, \yb, \pb)$ is a  $\Q$-rational function,
  explicitly known and the symmetry pseudo--group of $E_\fbar$ is zero-dimensional. The reader
  can find  in the appendix   \ref{annexe:diffalg} a brief introduction to  differential algebra. 

\subsection{The self-equivalence problem}

The system \pref{PDEsystem} is fundamental and 
can be treated by two different approaches : brute-force method based  on 
differential algebra (this section) and geometric approach relying  
on Cartan's theory of exterior  differential systems (the next section). 
While it is classically known  that the existence of at least one transformation ${\varphi:X \rightarrow X}$
 can be checked by  computing  the \emph{integrability conditions} of the system
\pref{PDEsystem}, which is completely algorithmic whenever the functions $f, \fbar:X\rightarrow \C$ are 
explicitly known \citep{BLOP}, there is no general algorithm for computing closed form of $\varphi$.
We shall show that if the symmetry pseudo--group of $E_\fbar$ is zero-dimensional, the transformation
 $\varphi$ is obtained \emph{without} integrating any differential equations.

\begin{defn}[Symmetry pseudo--group]
	To any differential equation $E_\fbar$ and any $\D$-groupoid~$\Phi$
	 that acts on $(x,y)$,  we associate the $\D$-groupoid 
	$\S_{E_\fbar,\Phi}$ formed by the formal Taylor series 
	solutions  of the \emph{self--equivalence} problem
	\begin{equation} \label{S:eqns}
         	E_\fbar = \sigma_*(E_\fbar) \mbox{ and } \sigma \in
		\Gamma \Phi.
	\end{equation}
	The symmetry pseudo-group $S_{E_\fbar, \Phi} := \Gamma \S_{E_\fbar, \Phi}$ is the set  of $C^\infty$-functions $\sigma: X 	
	\rightarrow X$ that are local solutions 
	of the Lie  defining equations \pref{S:eqns}.
\end{defn}

\begin{example} 
	Consider the $\D$-groupoid $\Phi := \Phi_3$ from example 
	\ref{ex} and the Emden-Fowler equation  $E_\fbar$ 
	(number 6. 11, page 544 in \citep{Kamke}) 
	\begin{equation}\label{emf:eq}
		y''= \Frac{1}{x y^2}.
	\end{equation}
	The Lie  defining equations of the symmetry pseudo-group of $E_\fbar$  	
	are obtained by setting  $f(x,y,p)=\Frac{1}{x y^2}$ and 
	$\fbar(\xb,\yb,\pb)=\Frac{1}{\xb \yb^2}$ in the equations  \pref{ex:f}.
	The characteristic set (see definition \ref{defn:cara} page 	
	\pageref{defn:cara}) of these equations is
	\begin{equation}\label{symemf}
		C_\sigma = 
		\left\{ \pb = p\ \Frac{\yb}{y},\ {\yb^3}=y^3,\ \xb = x \right\}.
	\end{equation}
	This PDE's system is particular :  it contains only non--differential 	equations.  We have 
	$\dim C_\sigma = 0$ and $\deg C_\sigma = 3$.
	The symmetry pseudo-group 
	$S_\fbar := S_{E_\fbar, \Phi}$ is actually a group with 3 elements:
	$$
		S_\fbar = \left\{ (x, y, p) \to (x, \lambda y, \lambda p)
			\mid \lambda^3 = 1 \right\}
	$$
	
\end{example}

\subsection{Equivalence problem with fixed target}

Assume that  $\fbar(\xb, \yb, \pb)$ is 	a fixed $\Q$-rational function. 

\begin{example}
   Consider again  the equivalence problem of example \ref{ex} :
   \begin{equation} \label{Sigmaexple}
	\begin{gathered}
	\bar p - \bar y_x - p\bar y_y=0,\\
	\bar y_{xx} + 2p\bar y_{xy} + p^2\bar y_{yy} + f\ \bar y_y 
	-\fbar(\bar x,\,\bar y,\,\bar p) =0,  \\
	\bar x_x-1=0,\, \bar x_y=0,\, \bar x_p=0,\, \bar y_p=0,\, \bar y_y\neq 0.
	\end{gathered}
   \end{equation} 
   These equations constitute a quasi-linear characteristic set w.r.t. the elimination
   ranking 
   $$
   	\Theta f \succ \Theta \pb \succ \Theta \yb \succ \Theta \xb.
   $$ 
   Consequently, the associated differential ideal is prime.
   This fact can be generalized to any  $\D$-groupoid $\Phi$ defined by quasi-linear characteristic set   (see proposition \ref{prop:prolongation}).
\end{example}

\begin{proposition} \label{cor:cible}	
	The PDE's system \pref{PDEsystem} (where  $\fbar(\xb, \yb, \pb)$ is
	a fixed rational function)
	is a quasi--linear characteristic set $C$
	w.r.t the  ranking $\Theta f \succ \Theta \pb \succ \Theta\{\bar y, \bar x\}$.
\end{proposition}

\subsection{Brute-force method} \label{brute:force}

Using  \RG\ we compute  
a new characteristic set $C$ of the PDE's system \pref{PDEsystem} w.r.t. the new elimination ranking
 $\Theta\{ \pb,\,\yb,\,\xb \} \succ \Theta\{ f \}$, which eliminates the indeterminates $\{ \pb,\,\yb,\,\xb \}$.
We make the partition of $C$ as in the formula \pref{def:partition}
\begin{equation} \label{bf:partition}
	C = C_f \sqcup C_\varphi
\end{equation}
where $C_f := C \, \cap\, \Q[x,y,p] \{ f \}$ and $C_\varphi := C \setminus C_f$.


\begin{proposition}
	The transformation $\varphi$ does exist for \emph{almost any} function $f$ 
	satisfying the PDE's system associated to the characteristic set $C_f$. 
	The function $\bar \x=\varphi(\x)$ is solution of the
	PDE's system associated to  $C_\varphi$. 
\end{proposition}

\begin{defn}
	When $\dim C_\varphi = 0$, the  algebraic system associated to 
	$C_\varphi$ is called 
	the \emph{necessary form of the change of coordinates} $\xbar = \varphi(\x)$.
\end{defn}

\begin{example}
Consider the equivalence problem of example \ref{ex}. Suppose that the target~$E_\fbar$ is the Airy 
equation 
$$
	\bar y''= \bar x \bar y.
$$
In this case, \RG\ returns $C_\varphi$ and $C_f$ resp. given by \pref{thairyetac} and  \pref{thairyf}
\begin{meqnarray}\label{thairyetac}
	\bar y_{{xx}} &=& -f \bar y_{{y}}+pf_{{p}}\bar y_{{y}}-\frac{1}{2}{p}^{2}f_{pp}
		\bar y_{{y}}+\bar y f_{{y}}-\frac{1}{2}\,\bar y f_{{xp}} -\frac{1}{2}\bar y f_{{pp}}f +\frac{1}{4}\bar y {f_{{p}}}^{2}-\frac{1}{2}\bar y pf_{{yp}}\\
	\bar y_{{xy}} &=& -\frac{1}{2}f_{{p}}\bar y_{{y}}+\frac{1}{2}pf_{{p,p}}\bar y_{{y}}\\
	\bar y_{{yy}} &=& -\frac{1}{2}f_{{pp}}\bar y_{{y}},\\
	\bar y_p &=& 0,\\
	\bar x &=& f_{{y}}-\frac{1}{2}f_{{xp}}-\frac{1}{2}f_{{pp}}f +\frac{1}{4}{f_{{p}}}^{2}-\frac{1}{2}\,pf_{{yp}}
\end{meqnarray}
\begin{meqnarray} \label{thairyf}
	f_{{xxp}} &=& 2\,f_{{xy}}+f_{{p}}f_{{xp}}
		-2+{p}^{2}f_{{yyp}}-f_{{pp}}f_{{x}} +\cdots\\ 
	f_{{xyp}} &=& 2\,f_{{yy}}-pf_{{yyp}}-f_{{ypp}}f -f_{{pp}}f_{{y}}+f_{{p}}f_{{yp}}\\
	f_{{xpp}} &=& f_{{yp}}-pf_{{ypp}}\\
	f_{ppp}   &=& 0.
\end{meqnarray}%
We have $\dim C_\varphi = 3$ which means that the transformation 
$\xbar = \varphi(\x)$, when $f$ satisfies~$C_f$, depends on 3 arbitrary constants.
\end{example}

\begin{example}
 Assume now that  
 the target equation  $E_\fbar$ is
$$
	 \bar y''= \bar y^3.
$$
 \RG\ returns  $C_\varphi$  and $C_f$  resp.  given by   \pref{eq:ex3} and \pref{f:ex3}
\begin{meqnarray} \label{eq:ex3}
	\bar y^{2} &=& \frac{1}{12}\,({4\,f_{{y}}-2\,f_{{xp}}-2\,f_{{pp}}f_		{{}}-2\,pf_{{yp}} +{f_{{p}}}^{2}}),\\
	\bar x &=& x,
\end{meqnarray}%
\begin{meqnarray}\label{f:ex3}
	f_{{xxxp}} &=&\left({4\,f_{{y}}-2\,f_{{xp}}-2\,f_{{pp}}f_{{}}-2\,pf_{{yp}}
		+{f_{{p}}}^{2}}\right)^{-1}\times\\
		 &&\left(24\,{p}^{2}{f_{{yp}}}^{2}f_{{y}}-24\,{p}^{2}f_{{yy}}f_{{p}}f_{{yp}}+ 
		  \cdots +12\,{p}^{2}{f_{{x,yp}}}^{2}+12\,{f_{{}}}^{2}{f_{{yp}}}^{2}-8\,{f_{{pp}}}^{3} {f_{{}}}^{3}\right)\\
	f_{{xxyp}} &=& ({4\,f_{{y}}-2\,f_{{xp}}-2\,f_{{pp}}f_{{}}-2\,pf_{{yp}}
		+{f_{{p}}}^{2}})^{-1} \times\\
		&& (-4\,pf_{{ypp}}f_{{}}f_{{p}}f_{{yp}}-4\,{f_{{xp}}}^{2}f_{{yp}}+  \cdots   +6\,pf_{{p}}f_{{y,yp}}f_{{pp}}f_{{}}+2\,{p}^{3}{f_{{yyp}}}^{2})\\
	f_{{xyyp}} &=& \,\,({4\,f_{{y}}-2\,f_{{xp}}-2\,f_{{pp}}f_{{}}-2\,pf_{{yp}}
		+{f_{{p}}}^{2}})^{-1} \times\\
		&&  (2\,{f_{{ypp}}}^{2}{f_{{}}}^{2}-2\,f_{{pp}}f_{{ypp}}f_{{}}pf_{{yp}}+  \cdots  +4\,{f_{{yp}}}^{2}f_{{xp}}+4\,{f_{{yp}}}^{3}p+16\,f_{{ypp}}{f_{{y}}}^{2}) \\
	f_{{xpp}} &=& f_{{yp}}-pf_{{ypp}}\\
	f_{{ppp}} &=& 0.
\end{meqnarray}%
Consequently $\dim C_\varphi = 0$ and $\deg C_\varphi = 2$.
Thus, $\varphi$ is an algebraic transformation of degree 2, given by equations~\pref{eq:ex3}.
\end{example}

\subsection{From equivalence problem with determined target $E_\fbar$ to the self-equivalence problem} \label{self:EPB}
Consider the characteristic set   $C = C_f \sqcup C_\varphi$ associated to the equivalence problem
 with determined target $E_\fbar$ defined by \pref{bf:partition} and computed w.r.t the 
  elimination ranking $\Theta\{ \pb,\,\yb,\,\xb \} \succ \Theta\{ f \}$.  On the other hand,   consider
  the characteristic set  $C_\sigma$ associated to the self-equivalence problem  \pref{S:eqns},
 computed w.r.t  the orderly ranking on  $\{ \pb,\,\yb,\,\xb \}$. By definition, we have
\begin{eqnarray*}
	C &\subset& \Q(x,y,p)\{ \xb,\,\yb,\,\pb, f \} \\
	C_\sigma &\subset& \Q(x,y,p)\{ \xb,\,\yb,\,\pb \}
\end{eqnarray*}

We obtain  $C_\sigma$ from $C$ by requiring that the two functions $f$ and  $\fbar$ are equal.  
The set $C$ is   \emph{specialized} by substituting the symbol $f$ by the value 
$\fbar(x,\, y,\, p)$.
After specialization, the differential system $C_f$ constraining the function $f$
is automatically satisfied since there exists at least one solution
$\xbar = \sigma(\x)$ of the problem, namely  $\sigma=\Id$. 

\begin{lemma} \label{lemme:base}
   The two characteristic sets $C_\varphi$ and $C_\sigma$ have  the    same dimension and
the same degree in the zero-dimensional case.
\end{lemma}

\begin{pf}
For each equation  $E_f$ equivalent to the target equation $E_\fbar$, denote by  $\Phi_{f,\fbar}$ 
the diffiety defined by  $C_\varphi$. The  $\D$-groupoid $\S_{E_\fbar,\Phi}$ acts simply transitively
 (see Figure  \ref{groupoidefig}) on the diffiety $\Phi_{f,\fbar}$, i.e.
$$
	\forall \varphi_0, \varphi \in \Gamma \Phi_{f,\fbar},
	\quad \exists ! \sigma \in \Gamma \S_{E_\fbar,\Phi},\quad \varphi = \sigma\circ \varphi_0.
$$
\begin{figure}[h]
	\begin{center}
\setlength{\unitlength}{0.00027333in}%
\begingroup\makeatletter\ifx\SetFigFont\undefined%
\gdef\SetFigFont#1#2#3#4#5{%
  \reset@font\fontsize{#1}{#2pt}%
  \fontfamily{#3}\fontseries{#4}\fontshape{#5}%
  \selectfont}%
\fi\endgroup%
{\renewcommand{\dashlinestretch}{30}
\begin{picture}(3311,2233)(0,-10)
\drawline(421,1672)(2349,270)
\drawline(2281.998,297.054)(2349.000,270.000)(2302.612,325.401)
\drawline(2775,1774)(2775,289)
\drawline(2755.845,365.610)(2775.000,289.000)(2794.155,365.610)
\drawline(582,1924)(2491,1924)
\drawline(2419.960,1906.240)(2491.000,1924.000)(2419.960,1941.760)
\put(2700,1849){\makebox(0,0)[lb]{{\SetFigFont{12}{14.4}{\rmdefault}{\mddefault}{\updefault}
\small{$\jet{\xbar_0}{} \fbar$}
}}}
\put(-950,1849){\makebox(0,0)[lb]{{\SetFigFont{12}{14.4}{\rmdefault}{\mddefault}{\updefault}
\small{$\jet{\x}{} f$}
}}}
\put(1425,2074){\makebox(0,0)[lb]{{\SetFigFont{12}{14.4}{\rmdefault}{\mddefault}{\updefault}
\small{$\varphi_0$}
}}}
\put(2850,874){\makebox(0,0)[lb]{{\SetFigFont{12}{
14.4}{\rmdefault}{\mddefault}{\updefault}
\small{$\sigma$}
}}}
\put(595,874){\makebox(0,0)[lb]{{\SetFigFont{12}{14.4}{\rmdefault}{\mddefault}{\updefault}
\small{$\varphi$}
}}}
\put(2475,-240){\makebox(0,0)[lb]{{\SetFigFont{12}{14.4}{\rmdefault}{\mddefault}{\updefault}
\small{$\jet{\xbar}{} \fbar$}
}}}
\end{picture}
}
	\caption {Simply transitive action of $\S_{E_\fbar,\Phi}$ on $\Phi_{f,\fbar}$
        where $\xbar_0 = \varphi_0(\x)$ and $\xbar=\varphi(\x)$}
	\label{groupoidefig}
	\end{center}
\end{figure}

Every $\varphi_0 \in \Gamma \Phi_{f,\fbar}$, define a bijective correspondence $\S_{E_\fbar,\Phi} \to \Phi_{f,\fbar}$
$$
   		\jet{\xbar_0}{\infty}\sigma \longrightarrow
   		\jet{\x}{\infty}\varphi = 
   			(\jet{\xbar_0}{\infty}\sigma)\circ (\jet{\x}{\infty}\varphi_0),
   			\quad (\sigma \in S_{E_\fbar,\Phi}).
$$   	
In fact, according to the Taylor series composition formulae, the  one-to-one correspondence between
 the two algebraic varieties $\S_{E_\fbar,\Phi}$ and $\Phi_{f,\fbar}$ is bi-rational. Consequently, these two varieties have the  same dimension  and
 the same degree in the zero-dimensional case. The same goes  for the two characteristic sets  $C_\varphi$ and $C_\sigma$ defining
 these varieties. 
\end{pf}	

We are now  in position to announce  the main theorem of the paper (recall that $\fbar$ is assumed $\Q$-rational
 function of its arguments and by definition, the degree of an algebraic transformation  $\bar \x=\varphi(\x)$ is the generic number
 of points $\xbar$ when $\x=(x,y,p)\in \C^3$ is determined):
 
\begin{thm}\label{thm2}
	The following conditions are equivalent 
	\begin{marray}
		(1) & $\dim(C_\varphi)=0$, \\
		(2) & $\dim(C_\sigma) =0$, \\
		(3) & $\dim(S_{E_\fbar,\Phi})  =0$, \\
		(4) & $\card(S_{E_\fbar,\Phi}) < \infty$.
	\end{marray}%
	In this case, $\card S_{E_\fbar,\Phi} = \deg(C_\varphi) = \deg \varphi$.
\end{thm}

\begin{remark}
When the transformation $\xbar := \varphi(\x)$ is locally bijective, but not
 globally,~$S_{E_f,\Phi}$ and $S_{E_\fbar,\Phi}$ need not to have the same degree. Indeed, consider again the
groupoid~$\Phi_3$ and the equations 
$$
	y''= \Frac{6y^4+x-2{y'}^2}{2y}  \mbox{ and } \bar y'' = 6\bar y^2 + \bar x
$$
which are equivalent under  $(\bar x=x, \  \bar y=y^2)$. The corresponding 
symmetry group are respectively given by 
$$
	S_{E_f,\Phi}= \{ (x,y) \to (x, \lambda y) \mid \lambda^2=1 \} \mbox{ and } 
	S_{E_\fbar,\Phi} = \{ \Id \}.
$$
They have the same dimension but different cardinality. 
\end{remark}
		
\subsection{Expression swell}

In practice,  the above brute--force method, which consists 
of applying \RG\, to the PDE's system \pref{PDEsystem},  is rarely
effective due to expressions swell. Most of the examples treated 
here and in \citep{dridi:these}, using our  algorithm {\sf ChgtCoords},
cannot be treated with this approach.

It seems that the problem lies in the fact that we can not separate
the computation of~$C_\varphi$ from that of $C_f$ which contains,
very often, long expressions (observe that, since $C_\varphi$ is computed
in terms of $f$ then there is no need to compute the equivalence conditions~$C_f$).

An other disadvantage of the above method is that  we have to restart computation from the very beginning
if  the target equation is changed. 
In the  next section, we propose our algorithm {\sf ChgtCoords} 
to compute the transformation $\varphi$ alone and in terms of differential
invariants. These invariants are provided by Cartan method for a generic 
$f$ which means that we have not to re-apply Cartan method if the target equation is changed and a big part of calculations is generic. 
Furthermore, the computation of $\varphi$ in terms of differential invariants significantly  reduces the size of the expressions.

\section{Cartan's method based approach}\label{cartan}

 We refer the reader to \citep{cartan:pb, Hsu2, Olver1,  neut:these, dridi:these}  for an
   expanded tutorial presentation on Cartan's equivalence method
  and application to second order ODE. 

When applied, Cartan's  method  furnishes a finite set of fundamental invariants
  and  a certain number of invariant derivations generating the field of the differential invariants.  

\begin{example}
Consider the equivalence problem of example \ref{ex}. The
PDE's system \pref{PDEsystem} reads ($p=y'$)
$$
\underbrace{
\begin{pmatrix}
\d\bar p-\fbar(\bar x,\bar y,\bar p)\d\bar x\\
\d\bar y-\bar p\d\bar x\\
\d\bar x
\end{pmatrix}}_{\omega_\fbar}
=
\underbrace{
\begin{pmatrix}
a_1 & a_2 & 0\\
0 & a_3 & 0\\
0 & 0 & 1\\ 
\end{pmatrix}}_{S(\a)}
\underbrace{
\begin{pmatrix}
\d p- f(x, y,p)\d x\\
\d y- p\d x\\
\d x
\end{pmatrix}}_{\omega_f}
$$
with $\det(S(\a))\neq 0$. In accordance with Cartan, this system is lifted to the new linear Pfaffian system 
$$
	S(\bar \a)\ \omega_\fbar=S(\a)\ \omega_f.
$$ 
After two normalizations and one prolongation, Cartan's method yields (an $e$-structure with)
 three fundamental invariants
\begin{equation} \label{invs}
\begin{array}{ll}
  I_1 = -\Frac{1}{4}(f_{p})^{2} - f_{y} + \Frac{1}{2}D_{x}f_{p}, \quad
  I_2 = \Frac{f_{ppp}}{2 {a}^{2}},  \quad
  I_3 = \Frac{f_{yp} - D_{x}f_{pp}}{2a},
\end{array}
\end{equation}%
and the invariant derivations
\begin{equation}\label{derivations}
\begin{array}{l}
	X_1 = \Frac{1}{a}\Frac{\partial}{\partial p}, \quad
	X_3 = D_x- \Frac{1}{2}f_{p}a\Frac{\partial}{\partial a}, \quad
	X_4 = a\Frac{\partial}{\partial a},\\[3mm]
	X_2 =  \Frac{1}{a}\Frac{\partial}{\partial y} 
    		+ \Frac{1}{2}\Frac{f_{p}}{a}\Frac{\partial}{\partial p} 
    		- \Frac{1}{2}f_{pp}\Frac{\partial}{\partial a},
\end{array}
\end{equation}
where  $a=a_3$ and as usual  $D_x=\frac{\partial }{\partial x} +  p\frac{\partial }{\partial y} 
          + f(x,y,p)\frac{\partial }{\partial p}$ denotes the Cartan vector field. 

When $\dim(\S_{E_\fbar,\Phi})=0$, the additional parameters $a$ and $\bar a$
can be (post)normalized by fixing some invariant to some suitable value.
\end{example}

\begin{proposition}[\cite{Olver1}]\label{tmOlver}
	The symmetry groupoid  $\S_{\fbar, \Phi}$ is 
	zero-dimensional 	if and only if there exist exactly 
	three functionally independent specialized invariants
	$I_1[\fbar], I_2[\fbar], I_3[\fbar]:X \rightarrow \C$. 
\end{proposition}

Note that the  invariants  $I_1[\fbar], \cdots, I_3[\fbar]$ are functionally
independent if and only if  $\d I_1[\fbar]\wedge \cdots\wedge
\d I_3[\fbar] \neq 0$ and  if  the function $\fbar$ is rational,  then 
the specialized invariants $I[\fbar]: X \to \C$ are algebraic functions in  $(x,y,p)$.

\subsection{Computation of $\varphi$}

Suppose that the symmetry groupoid $\S_{E_\fbar,\Phi}$ is zero-dimensional.
 Then, according to  theorem \ref{tmOlver}, there exists $3$  invariants
$F_k := I_k[\fbar],$  $1 \leq k \leq 3$ such that  the 
 algebraic (non differential) system
\begin{equation}\label{sys:direct}
    \left\{ F_1(\xbar) =  I_1, F_2(\xbar) =  I_2,  F_{3}(\xbar) = I_3 \right\}
\end{equation}
is locally invertible and has a finite number of solutions 
\begin{equation} \label{sys:ivrs}
	\xbar = F^{-1}(I_1,\ldots, I_3). 
\end{equation}
The specialization of $I_1, \ldots , I_3$ on the source function  $f$ yields
\begin{equation}\label{cartan:phi}
	\bar \x= F^{-1}(I_1[f],\,\ldots, I_3[f]).
\end{equation}

The main idea here is that the inversion \pref{sys:ivrs} is done by 
computing a (non differential) characteristic set 
$C$ for  the 
 system \pref{sys:direct} w.r.t. the ranking 
$\{ \bar x,\, \bar y,\, \bar p \} \succ \{I_1,\, I_2,\,I_3 \}$.  Now, 
the most simple situation happens when $\deg(C)=1$.
In this case, the necessary form of the change of coordinates 
$\varphi$ is the rational transformation defined by~$C$. 

\begin{example}
Consider the equivalence problem of the example  \ref{ex}  and the target equation~$E_\fbar$
introduced by G. Reid \citep{Reid93}
$$
	\bar y''=\Frac{\bar y'}{\bar x} + \Frac{4\bar y^2}{\bar x^3}.
$$  
The following invariants are functionally independent (where
$I_{i;j\cdots k} =X_k\cdots X_j(I_i)$ )

$$
\begin{gathered}
 \bar I_{1}=\Frac{3}{4\bar x^2} +8\,{\Frac {\bar y}{{\bar x}^{3}}},\quad
 \bar I_{1;3} ={\Frac {-3\,\bar x-48\,\bar y+16\,\bar p\bar x}{2{\bar x}^{4}}},\quad
 \bar I_{1;23} =-20{\Frac {1}{\bar a{\bar x}^{4}}},\quad
 \bar I_{1;31}=8{\Frac {1}{\bar a{\bar x}^{3}}}.
\end{gathered}
$$
We normalize the parameter $\bar a$ by setting $\bar I_{1,23}=-20$. 
The characteristic set $C$ is
\begin{system}\nonumber
\bar p  =-{\Frac {3}{32}}+{\Frac {3}{512}} {I_{1;31  }}^{2}I_1+{\Frac {1}{4096}} 
I_{1;3  }{I_{1;31  }}^{3},\\[5mm]\bar y  =-{\Frac {3}{256}} I_{1;31  }+{\Frac 
{1}{4096}} I_{I_{1;31 }}^{3},\\[5mm]\bar x  =\Frac{1}{8} I_{1;31 },
\end{system}%
which gives the sought necessary form of $\varphi$. As a byproduct we deduce that the symmetry
group $S_{E_\fbar,\Phi} = \{ \Id \}$. 
\end{example}

When  $\deg(C)$ is strictly bigger than 1, we have two cases.  
First, $\deg(C)=\deg(S_{E_\fbar,\Phi})$ and then $\varphi$ is the algebraic transformation defined
by $C$. Second,    $\deg(C) > \deg(S_{E_\fbar,\Phi})$. 
In this  case, to obtain the transformation~$\varphi$, we have to look for~3 other
functionally independent  invariants such that the new characteristic set~$C$ has
degree equal to $\deg(S_{E_\fbar,\Phi})$.

\begin{example}
Consider the equivalence problem of example \ref{ex} and the target equation~number 6.9 in  \citep{Kamke}:
$$
	\bar y''=\bar y^3 + \bar x\bar y.
$$
The corresponding  symmetry group  is 
$$
	S_{E_\fbar,\Phi} =\left\{(x,y) \into (x,\,\lambda y) \, |
          \lambda^2=1 \right\}.
$$
One can verify that  $I_1,\, I_{1;13}$ and $I_{1;133}$,
when specialized on the considered equation, are functionally independent. 
The associated  characteristic set $C$~is
\begin{system}\nonumber
\bar p &=&- {\Frac { \left(4 {\bar x}^{2}+ 2 {I_1}\bar x-3 {I_{1;33}}
        -2 {I_1}^{2} \right)}{ 3(I_{1;3}+1)}} \bar y,\\
{\bar y}^{2}&=& - \Frac{1}{3}\bar x-\Frac{1}{3} {I_1},\\[3mm]
{\bar x}^{3} &=&  -\Frac{3}{2}{I_1} {\bar x }^{2}  +\Frac{3}{4} {I_{1;33}} \bar x 
              - \Frac{3}{4}{ I_{1;3}}-\Frac{3}{8} {{ I_{1;3}}}^{2}+ \Frac{3}{4}{I_{1;33}} {I_1} + \Frac{1}{2} {{I_1}}^{3}-\Frac{3}{8},
\end{system}%
which has degree 6, strictly bigger than the degree of the symmetry group.
However, if  we consider  the invariants  
$K_1 := {I_{1;233}}/{I_{1;31}}$, $K_2 :={I_{1;234}}/{I_{1;31}}$ and
$K_3 := {I_{1:231}}/{I^2_{1;31}}$, we obtain 
\begin{system} \nonumber
	\bar p &=& -K_1\bar y,\\
	{\bar y}^{2} &=& \Frac{1}{6}K_3,\\
	\bar x &=& - \Frac{1}{6}K_3 +K_1.
\end{system}%
That is,  the necessary form of $\varphi$ since  this new set  has  degree two.

\end{example}

\subsection{Heuristic of degree reduction}
We have to remark here that searching invariants  giving the required degree (as in example above)
 is not  an easy task, although algorithmic. This is simply because  the algebra of invariants can be very large.
For this reason we provide an  important heuristic which enables us to obtain the desired degree
even for a "bad" choice of invariants.

\begin{example}
Consider the Emden-Fowler equation \pref{emf:eq} and the $\D$-groupoid of
transformations $\Phi_3$. We have already computed the corresponding symmetry
groupoid. The specialization of the invariants  $I_1,\, I_{1;13}$ and $I_{1;133}$ gives three 
functionally independent functions. As explained
above,  we obtain the following characteristic set computed w.r.t. the ranking
$\bar p \succ \bar y \succ \bar x \succ  I_1 \succ  I_{1;3} \succ
I_{1;33}$

\begin{system}\label{3resG1}
\bar p &=& \left( \Frac{3}{8}I_{1}-
\Frac{1}{4}{\Frac {I_{1;33}}{I_{1}}}+ \Frac{1}{3} {\Frac
{{I_{1;3}}^{2}}{{I_{1}}^{2}}} \right)\bar x\bar y - \Frac{1}{6} {\Frac 
{I_{1;3}}{I_{1}}}\bar y,  \\\bar y^3 &=& \left( -\Frac{9}{4}-2 {\Frac {{{\it 
I_{13}}}^{2}}{{I_{1}}^{3}}}+ \Frac{3}{2} {\Frac {I_{1;33}}{{I_{1}}^{2}}} \right) 
\bar x-{\Frac {I_{1;3}}{{I_{1}}^{2}}},   \\\bar x^2 &=& 4 \left( {\Frac {I_{1;3} 
I_{1} }{9 {{\it I_1}}^{3}-8 {I_{1;3}}^{2}+6 I_{1;33} I_{1}}} \right)\bar x  
+ 8 {\Frac {{I_{1}}^{2}}{9 {I_{1}}^{3}-8 {I_{1;3}}^{2}+6 I_{1;33} I_{1}}}. 
\end{system}%

Comparing with the $\D$-groupoid of symmetries \pref{symemf} we deduce that, in
contrast to $\bar y$, the degree of $\bar x$ must be reduced to one.
This can be done  in the following manner. First, observe that the Lie defining
equations of $\Phi_3$,  more exactly $\bar x_p=0$, implies that~${X_1(\bar x)=0}$
where $X_1=\frac{\partial}{a\partial p}$ is the invariant derivation~\pref{derivations}. Now,
differentiate the last equation of the  characteristic set, which we write as
$\bar x^2=A\bar x+B$, w.r.t the derivation $X_1$. We find $A_{;1}\bar x +
B_{;1}=0$.
The coefficient of $\bar x$ in this equation, which is
invariant, could not vanish (since  it is  not identically zero when specializing on
the Emden--Fowler equation). We obtain  
$\bar x = - \Frac{  B_{;1} }{ A_{;1}  }$
or explicitly
\begin{equation}\label{xbar}
	\bar x= -2 {\Frac { K I_{1;1}  + I_1K_{;1}  }{ K I_{1;31} +{\it I_{1;3}}
  	K_{;1}  }}\ \mbox{\ with\ } K={\Frac {{{\it I_1}}}{9 {{\it I_1}}^{3}-8 {{\it
	I_{1;3}}}^{2}+6 {\it I_{1;33}} {\it I_1}}}.
\end{equation}
The necessary form of
the change of coordinates $\varphi$ is then given  by 
 \pref{xbar} and 
the  first two equations of
\pref{3resG1}.
\end{example}
The above reasoning can be summarized as follows 

\begin{center}
\vspace{10pt}
\begin{tabular}{|p{0.95\linewidth}|}
  \hline
   {\sc Procedure} {\sf ChgtCoords} \\
   {\bf Input :} $E_\fbar$ and $\Phi$ such that $\dim(S_{E_\fbar,\Phi})=0$ \\
   {\bf Output :} $\bar \x=\varphi(\x)$ the necessary form of the change of coordinates \\
  \hline
   1- Find $3$ functionally independent invariants $(I_1[\fbar],\, I_2[\fbar],  \, I_3[\fbar])$
    defined on $X$.\\
   2- Compute a char. set $C$ of the algebraic system \pref{sys:ivrs}. \\
   3- If $\deg(C)=1$ then Return $C$. \\
   4- Compute $S_{E_\fbar,\Phi}$ with \RG.\\
   5- WHILE $\deg(C) \neq \deg(S_{E_\fbar,\Phi})$ DO \\
    ~~~~ Reduce the degree of $C$.  \\
    END DO\\
   6- Return $C$.\\
  \hline
\end{tabular}
\vspace{10pt}
\end{center}

\section{The solver}

\subsection{Pre-calculation of $\varphi$}
\subsubsection{The first step : the adapted $\D$-groupoid}

Let $\Phi_1,\,\ldots,\, \Phi_7 \subset \J^\infty_{*}(\C^2,\C^2)$ denote  the $\D$-groupoids
 defined in  the 
table \ref{table} page \pageref{table}. It is not difficult to see that 
$\Phi_1\subset\Phi_3\subset\Phi_5$ and  
$\Phi_2\subset\Phi_4\subset\Phi_6$ and finally 
$\Phi_5, \Phi_6 \subset \Phi_7$.


\begin{defn}[Signature index]
	The \emph{signature index} of  $E_f$ is
	$$
	\sign(E_f) :=\left( (\d_1, \d_3, \d_5),(\d_2, \d_4,\d_6), \d_7\right)
	\mbox{ where } \d_i := \dim S_{E_f,\Phi_i},\ 1\leq i\leq 7.
	$$
\end{defn}
Clearly, $(\d_1 \leq \d_3 \leq \d_5 \leq \d_7)$ and $(\d_2 \leq \d_4
\leq \d_6 \leq \d_7)$.  Recall  that the calculation of
theses dimensions  does not require  solving differential equations.
\begin{defn}
	We shall say that the signature index $\mathrm{sign}(E_f)$ \emph{matches} the signature index
	$\sign(E_\fbar)$ if and only if
	$$
		\d_7=\bar \d_7 \mbox{ and }  (s_1=\bar s_1 \mbox{ or } s_2=\bar s_2)
	$$ where $s_1$ and $s_2$ stand for $(\d_1, \d_3, \d_5)$ and $(\d_2, \d_4,\d_6)$ resp. 
\end{defn}
\begin{defn}
Two second order ODE $E_f$ and $E_\fbar$ are said to be \emph{strongly equivalent} if 
$$
	\exists \Phi\in\{ \Phi_1,\cdots, \Phi_7\}, \
	\exists \varphi\in \Gamma\Phi,\ \varphi_*E_f=E_\fbar,\ \dim S_{E_\fbar,\Phi}=0.
$$
\end{defn}
\begin{lemma}
If $E_f$ and $E_\fbar$ are strongly equivalent then their signature indices
match. 
\end{lemma}
\begin{defn}[Adapted $\D$-groupoid] \label{def:adapted}
	A $\D$-groupoid $\Phi$ is said to be \emph{adapted} to the differential equation~$E_\fbar$ 
	if $\dim (S_{E_\fbar,\Phi})=0$ and~$\Phi$ is maximal among  $\Phi_1,\,\cdots,\, \Phi_7$ satisfying this property.
\end{defn}

\begin{center} \scriptsize
\begin{table}
\begin{tabular}{|c|l|p{8cm}|} \hline
&{Transformations}  & Equation number according to Kamke's book \\ \hline    
$\Phi_1$ & $ \bar  x=x,\ \bar y=\eta(x,y)$
	& 1, 2, 4, 7, 10, 21, 23, 24, 30, 31, 32, 40, 42, 43, 45, 47, 50 \\ \hline    
$\Phi_3$ & $  \bar  x=x+C,\ \bar y=\eta(x,y)$ 
	&11, 78, 79, 87, 90, 91, 92, 94, 97, 98, 105, 106, 156, 172 \\ \hline
$\Phi_5$ & $ \bar  x=\xi(x),\ \bar y=\eta(x,y)$ & Null \\ \hline
$\Phi_2$ & $ \bar  x=\xi(x,y),\ \bar y=y$ &  81, 89, 133, 134, 135, 237 \\ \hline
$\Phi_4$ & $ \bar  x=\xi(x,y),\ \bar y=y+C $  
	& 11, 79, 87, 90, 92, 93, 94, 97, 98, 99, 105, 106, 172, 178 \\ \hline
$\Phi_6$ & $ \bar  x=\xi(x,y),\ \bar y=\eta(y)$  
    & 80, 86, 156, 219, 233  \\ \hline
$\Phi_7$ & $ \bar  x=\xi(x,y),\ \bar y=\eta(x,y)$  
    & 3, 5, 6, 8, 9, 27,  44, 52,  85, 95, 108, 142, 144, 145, 147, 171, 211, 212, 238 \\ \hline
\end{tabular}\\
\caption{Adapted groupoids for certain equations from Kamke's list} \label{table}
\end{table}
\end{center}

The   table \ref{table} associates to each equation in the third column\footnote{A more complete list of equations is available upon request.}
 its  adapted groupoids.  For instance, the first Painlev\'e equation (number 3) appears 
in the last row which means that its adapted $\D$-groupoid is 
the point transformations $\D$-groupoid $\Phi_7$ . To the  Emden--Fowler  equation,  
number~11, we associate the $\D$-groupoids  $\Phi_3$  and $\Phi_4$. In the case 
of homogeneous linear second order ODE (e.g. Airy equation, Bessel equation, 
Gau\ss\ hyper-geometric equation) we prove that, generically, the 
adapted $\D$-groupoid is~$\Phi_4$.

\subsubsection{The second step}
Once the list of adapted $\D$-groupoids $\Phi$  is known, we proceed by computing the necessary form  of the change of 
coordinates $\varphi\in \Gamma\Phi$ using {\sf ChgtCoords}. Doing so, we construct a \Maple\ table indexed by  Kamke's book 
equations and where  entries corresponding to the index $E_\fbar$ are: 
\begin{marray}
 1- & the signature index of $E_\fbar$, \\
 2- & the list of the adapted $\D$-groupoids $\Phi$ of  $E_\fbar$,\\ 
 3- & the necessary form of the change of coordinates $\varphi\in \Gamma\Phi$.
\end{marray}
For instance, the entries associated to Rayleigh equation $y''+ {y'}^4 + y=0$ are: 
\begin{marray}
 1- & the signature index $\left((0,1,1),(1,1,1),1\right)$, \\
 2- & the $\D$-groupoid~$\Phi_1$, \\
 3- & the necessary form of the change of coordinates 
\end{marray} 
\begin{system}\nonumber
	\bar p &=& -36{\Frac {I_{2;1}}{72+72I_1+
	{I_{2;1}}^{2}}}\bar y,\\
	\bar x &=&x,\\
	{\bar y}^{3} &=& {\Frac {-1}{559872 {I_{2;1}^{2}}} }(I_{2;1}^{6}+216I_1I_{2;1}^{4}+ 216I_{2;1}^{4}  
	                    +15552I_1^{2}I_{2;1}^{2}+  31104I_1I_{2;1}^{2} +373248I_1^{3}\\
	          && 	+ 15552I_{2;1}^{2} +1119744I_1^{2} +1119744I_1 +373248 )
\end{system}%
with the normalization 
$
{I_2}/{I_{2;1}}=1.
$
Invariants here are those generated by  \pref{invs} and \pref{derivations} plus 
the essential invariant $\bar x=x$.

\subsection{Algorithmic scheme of the solver}

To integrate a differential equation $E_f$ our solver proceeds as follows
\begin{center}
\vspace{10pt}
\begin{tabular}{|p{0.95\linewidth}| }
  \hline
   {\sc Procedure} {\sf Newdsolve}  \\
   {\bf Input :} $E_f$  \\
   {\bf Output :} An equation $E_\fbar$ in  Kamke's list and the transformation 
   $\varphi$ such that $\varphi_*(E_f)=E_\fbar$ \\ \hline
   1- Compute the signature index of $E_f$.\\
   2- Select from the table the list of equations $E_\fbar$ such that $\sign(E_\fbar)$
   matches $\sign(E_f)$.\\
   3- FOR each equation $E_\fbar$ in the selected list DO \\
   	~(i) Specialize, on $E_f$, the necessary form of the change of coordinates 
   	associated to~$E_\fbar$. We obtain~$\varphi$.\\
   	~(ii) If $\varphi \in \Gamma\Phi$ and $\varphi_*(E_f)=E_\fbar$ 
     	then return ($E_\fbar$, $\varphi$). \\
    END DO. \\
  \hline
\end{tabular} 
\vspace{10pt}
\end{center}

\subsection{Features of the solver}
It is worth noticing  that the  time required to perform steps (i)- (ii) is very small. 
In fact, it is about one hundredth of a  second using  Pentium(4) with 256 Mo.
Experimented on many examples, the total time needed to solve a given equation does not
exceed few seconds in the worse situations.

The second feature of our solver is,  contrarily  to the symmetry methods, 
neither the table construction nor the algorithm of the solver involves integration of differential equations.  
Indeed, even the computation of signature indices is performed without solving the Lie equations.
\appendix

\section{Appendixes}

\subsection{Differential algebra} \label{annexe:diffalg}

The reader is assumed to be familiar with the basic notions  and 
notations of  differential algebra. Reference books are
\citep{Ritt} and \citep{kolchin:livre}. We also refer to \citep{BLOP, hubert00, Boulier}.
Let $U=\{u_1,\cdots,u_n\}$ be a set of differential 
indeterminates. $\kk$ is a differential field of 
characteristic zero endowed with the set of derivations 
$\Delta=\left\{ \partial_1,\,\cdots,\partial_p \right\}$. The monoid of derivations
\begin{equation}
	\Theta := \left\{ \partial_1^{\alpha_1} \partial_2^{\alpha_2} \cdots
		\partial_p^{\alpha_p} \mid \alpha_1, \ldots, \alpha_p \in \N \right\}
\end{equation}
acts freely on the alphabet $U$ and defines a new (infinite) alphabet $\Theta U$.
The differential ring of the polynomials built over $\Theta U$ with 
coefficients in $k$ is denoted $R=\kk\{U\}$.
Fix an admissible ranking over $\Theta U$. For $f \in R$, 
$\leader(f) \in \Theta U $ denotes the \emph{leader} (main variable), 
$I_f \in R$ denotes the \emph{initial} of $f$ and $S_f\in R$ denotes
the separant of $f$. Recall that $S_f=\frac{\partial f}{\partial v}$ where $ v = \leader(f)$. 
Let $C\subset R$ be a finite set of differential polynomials. 
Denote by $[C]$ the differential ideal generated by $C$ and by $\sqrt{[C]}$ 
the radical of $[C]$. Let $H_C := \{ I_f \mid f \in C\} \cup \{ S_f \mid f \in C \}$.
As usual, $\remf$ is the Ritt full reduction algorithm \citep{kolchin:livre}. 
If $r = \remf(f,C)$ then $\exists h \in H^\infty_C,\, hf=r \mod [C]$. 
Then the \emph{reduced form} is defined by $\RF(f) := r/h$.
\begin{defn}[Characteristic set]\label{defn:cara}
	The set $C\subset R$ is said to be a \emph{characteristic set} of the differential ideal
	${\mathfrak{c} := \sqrt{[C]}:H_C^\infty}$ if  
	\begin{marray}
		(1)  & $C$ is auto-reduced, \\
		(2)  & $f\in \mathfrak{c}$ if and only if $\remf(f,\, C) = 0$.
	\end{marray}
\end{defn}

\begin{defn}[Quasi--linear characteristic set]
	The  characteristic set $C\subset R$ is said to be
	\emph{quasi--linear} if for each $f\in C$
	we have $\deg(f,\,v) = 1$ where~$v$ is the leader of~$f$.
\end{defn}

\begin{proposition} \label{C:C0}
	When the characteristic set $C$ is quasi--linear, the
	differential ideal $\mathfrak{c} := \sqrt{[C]}:H_C^\infty \subset R$ is prime.	
\end{proposition}
	
\subsection{Taylor series solutions space}

Let $\kk :=\C(x_1,\cdots,x_p)$ be the differential field of coefficients endowed with the set
of derivations 
$\left\{ \frac{\partial}{\partial x_1}, \cdots, \frac{\partial}{\partial x_p} \right\}$. 
Let $C$ be a characteristic
set of a prime differential ideal $\mathfrak{c} \subset R$.
We associate to $C$ the  system
\begin{equation} \label{C:PDEsyst}
	(C=0, H_C\neq 0)
\end{equation} 
of equations $f=0, \, f \in C$ and inequations $h \neq 0,\, h \in H_C$.
\begin{defn} [Taylor series solution]\label{TS:solns}
	A \emph{Taylor series solution} (with coefficients in $\C$) 
	of the PDE's system \pref{C:PDEsyst}
	 is a morphism $\mu: R\into \C$
	of (non differential) $\C$-algebras  such that
	$$[C] \subset \ker \mu \mbox{ and } H_C \cap \ker \mu =\emptyset.$$
\end{defn}
The  source of the Taylor solution $\mu$ is
$\source(\mu) := \left( \mu(x_1),\ldots, \mu(x_p) \right) \in \C^p$ and the target is
$\target(\mu) := \left( \mu(u_1),\ldots, \mu(u_n) \right) \in \C^n$.
The diffiety associated to the characteristic set   $C$ is the set of the formal Taylor solutions
of the system  \pref{C:PDEsyst}. 

The dimension of the solutions space of \pref{C:PDEsyst} is
the number of arbitrary constants appearing in the Taylor series solutions $\mu$ when the
source point $\x := \source(\mu) \in \C^p$ is determined.
Let $K$ be the fractions field $\Fr (R/\mathfrak{c})$.
Recall that the \emph{transcendence degree} of a field extension $K/\kk$ is the greatest number
of elements in $K$ which are $\kk$-algebraically independent. 
The degree $[K:\kk]$ is the dimension of $K$ as a
$\kk$-vector space. When $\trdeg (K/\kk)=0$, the field $K$ is algebraic over $\kk$ and ${[K:\kk] < \infty}$. 
If $f \in C$, we denote $\rank(f) := (v,d)$ where $v := \leader f$ and $d := \deg(f, v)$.
Let
\begin{eqnarray*}
	\rank C   &:=& \left\{ \rank(f)\mid f \in C \right\} \\
	\leader C &:=& \left\{ \leader(f)\mid f \in C \right\} \\
	\dim C &:=& \card \left( \Theta U \setminus \Theta(\leader C) \right) \\
	\deg C &:=& \Prod{f\in C}{} \deg(f, \leader f).
\end{eqnarray*}
\begin{proposition} \label{C:solns}
	$\dim C = \trdeg(K/\kk)$ is the dimension of the solutions space of \pref{C:PDEsyst}.
	If $\dim C = 0$ then the cardinal of the solutions space is finite
	and equal to $\deg C = [K:\kk]$.
\end{proposition}

\subsection{Differential elimination} \label{sect:elim}

Let $U=U_1 \sqcup U_2$ be a partition of the alphabet $U$.
A ranking which eliminates the indeterminates of $U_2$ is such that
\begin{equation}
	\forall v_1 \in \Theta U_1,\, \forall v_2 \in \Theta U_2, \quad  v_2 \succ v_1.
\end{equation}

Assume that $C$ is a characteristic set of 
the prime differential ideal $\mathfrak{c} = \sqrt{[C]}:H_C^\infty$ w.r.t.
the elimination ranking $\Theta U_2 \succ \Theta U_1$.
Let $R_1 := k\{ U_1 \}$ be the differential polynomials $k$-algebra generated by the set
$U_1$. Consider the set $C_1 := C \cap R_1$ and
the differential ideal $\mathfrak {c}_1 := \mathfrak{c} \cap R_1$.

\begin{proposition}\label{elim:prop}	
 	$C_1$ is a characteristic set of $\mathfrak{c}_1$.
\end{proposition}

Consider the differential field of fractions  $K := \Fr (R/\mathfrak{c})$ and
denote by $\alpha: R \to K$ the canonical $k$-algebra morphism. 
Let $K_1$ be the differential subfield of $K$ generated by the set $\alpha(R_1)$. 
Then $K_1$ is the fraction field associated to the prime differential ideal 
$\mathfrak{c}_1 := \mathfrak{c} \cap R_1$.
The partition of the  characteristic set
\begin{equation} \label{def:partition}
	C = C_1 \sqcup C_2 \quad (\mbox{i.e. } C_2 := C \setminus C_1).
\end{equation}
enables us to study the field extension $K/K_1$. 
\begin{proposition}
   	$\trdeg (K/K_1) = \dim C_2$. 
	If $\dim C_2 =0$ then $[K:K_1] = \deg C_2$.
\end{proposition}

\subsection{Groupoids} \label{annexe:groupoids}
\begin{defn}[Groupoid]
	A \emph{groupoid} is a category in which every arrow is invertible.
\end{defn}
Let $(\Phi, X, \circ, \source, \target)$ be a category.
Each arrow $\varphi \in \Phi$ admits a source $\source(\varphi) \in X$ and a
target $\target(\varphi) \in X$  which are \emph{objects} of this category.
The composition   $\varphi_2 \circ \varphi_1$ of the two arrows   $\varphi_1$ and $\varphi_2$ 
are defined when  $\target(\varphi_1)= \source(\varphi_2)$. 

If  $\Phi$  is a groupoid, for each arrow $\varphi \in \Phi$,
there exists a unique inverse arrow $\varphi^{-1}$ such that 
$\varphi^{-1}\circ\varphi=\Id_{\source(\varphi)}$ and 
$\varphi \circ \varphi^{-1}=\Id_{\target(\varphi)}$.

Let $X$ and $U$ be two manifolds and $\x \in X$.
The Taylor series up to order $q$ (i.e. the jet of order $q$) of a
function  $f:X \to U$, of class $C^q$,  is denoted $\jet \x q f$. 
The Taylor series  of~$f$ about  $\x$ is denoted $\jet \x {} f$ or $\jet \x \infty f$.
We shall say that  $\x \in X$ is the source and  $f(\x)\in U$ is the target of the $q$-jet $\jet x q f$. 

\begin{example} 
	For instance, when  $X=U=\C$, we have
	$$
	\jet x q f := \left( x, f(x), f'(x), \ldots, f^{(q)}(x) \right) \in \C^{q+2}.
	$$
	This jet is said to be \emph{invertible}
	if $f'(x) \neq 0$. 
	The jet of the function $\Id$ about the  point~$ x$ is $(x,x,1,0,\ldots, 0)$.
\end{example} 

For each  integer $q \in \N$ and each $\x \in X$, we set 
$\J_\x^q (X,U) := \bigcup_f \jet \x q f$. We denote by  $\J^q (X,U) := \bigsqcup_{\x \in X}\J_\x^q (X,U)$ the jets space up to order $q$.
We denote by $\J_*^q(X,X)$ the submanifold of $\J^q(X,X)$ formed by the  invertible jets.
Recall that $\J_*^q(X,X)$ is a groupoid \citep{OlverGroupoid} for the composition of Taylor series up to order~$q$ according to 
\begin{equation} \label{def:comp}
	\jet{\x}{q} (g \circ f) = \left( \jet{f(\x)}{q} g \right) \circ 
			\left( \jet{\x}{q} f \right).
\end{equation}

By definition, a  $\D$-groupoid \citep{Malgrange} $\Phi \subset \J_*^\infty(X,X)$
is a sub--groupoid of $\J_*^\infty(X,X)$ formed by the Taylor series 
solutions (see def. \ref{TS:solns}) of an algebraic PDE's system called 
the \emph{Lie defining equations}. This  system 
contains an inequation which  expresses   the  invertibility of  the jets.  

The set of $C^\infty$-functions $\varphi:X\to X$ that are local solutions of the Lie
defining equations of $\Phi$ is a \emph{pseudo-group} denoted by $\Gamma \Phi$. We define
$\dim \Gamma \Phi = \dim \Phi/X := \dim C$ and, 
if $\dim C = 0$, $\card \Gamma \Phi = \card \Phi/X := \deg C$  where~$C$ is a characteristic
set (see sect. \ref{diffalg}) of the Lie defining equations of $\Phi$. 

\subsection{Prolongation algorithm} \label{prolong}

Our aim, here, is to give an efficient way to prolong the  action of $\Phi^{(0)}$ on $\J^0(\C, \C)$
on the jets space $\J^n(\C, \C)$. For each integer $q \geq 0$, define the differential field
$$
	\kk^{(q)} := \Q(x,\, y,\, y_1,\ldots,\, y_q)
$$
and the ring of differential polynomials
$$
	R^{(q)}   := \kk^{(q)}\{ \bar x,\, \bar y,\, \bar y_1, \ldots,\, \bar y_q \}
$$
The differential field $\kk^{(q)}$ is the coefficients field of $R^{(q)}$
endowed with the set of derivations $\left\{ \frac{\partial}{\partial x},\, 
\frac{\partial}{\partial y},\,\ldots,\, \frac{\partial}{\partial y_q} \right\}$.
Let us assume that the  Lie defining equations of $\Phi^{(0)}$ are given by 
a quasi-linear characteristic set $C^{(0)} \subset R^{(0)}$. 
The  $\D$-groupoid  $\Phi^{(q)}$ acting on  $\J^q(\C, \C)$ and  prolonging the action of $\Phi^{(0)}$
 is characterized by a characteristic set $C^{(q)} \subset R^{(q)}$.
The prolongation formulae \citep{Olver} of the  point transformation 
$(x,y) \to \left( \xi(x,y),\, \eta(x,y) \right)$ are of the form
$$
	\bar y_q = \eta_q(x,\,y,\ldots,\, y_q), 
$$
where $\bar y = \eta(x,\, y)$ if $q=0$.
The computation of the characteristic set  $C^{(q)}$ is done incrementally using 
the infinite Cartan field $D_x := \PD{}{x} + y_1 \PD{}{y} + y_2 \PD{}{y_1} + \cdots$
\begin{eqnarray*}
        \eta_{q} &:=& D_x \eta_{q-1} \cdot {(D_x \xi)}^{-1} \\
        C^{(q)}  &:=& C^{(q-1)} \cup \left\{ \bar y_{q} - 
        	\RF\left( \eta_{q},\, C^{(q-1)}\right) \right\}
\end{eqnarray*}

\begin{proposition} \label{prop:prolongation}
	If $C^{(0)}$ is a \emph{quasi-linear} characteristic set of $\Phi^{(0)}$
	then  $C^{(q)}$ is a \emph{quasi-linear} characteristic set of $\Phi^{(q)}$
	w.r.t. the elimination ranking 
	$$
		\Theta \bar y_q \succ \Theta \bar y_{q-1} \succ \cdots \succ \Theta\{\bar y, \bar x\}.
	$$
\end{proposition}
The previous proposition gives an efficient method to prolong a $\D$-groupoid $\Phi$ without 
explicit knowledge of (the form of) the transformations $\varphi \in \Gamma \Phi$. 

\begin{ack}
We are thankful to Rudolf Bkouche and François Boulier for many useful discussions  and to the referees  for their valuable comments and suggestions.
\end{ack}

\bibliographystyle{elsart-harv}
\bibliography{c}

\end{document}